\documentclass{amsart}

\usepackage{hyperref, color}

\newtheorem{theorem}{Theorem}[section]
\newtheorem{lemma}[theorem]{Lemma}

\theoremstyle{definition}

\newtheorem{proposition}[theorem]{Proposition}
\newtheorem{corollary}[theorem]{Corollary}

\theoremstyle{remark}
\newtheorem{remark}{\bf Remark}

\numberwithin{equation}{section}

\newcommand{\field}[1]{\mathbb{#1}}
\newcommand{\real}{\field{R}}
\newcommand{\sect}{\overline{Sect}_{rad}}

\newcommand{\al} {\alpha}

\newcommand{\de} {\delta}       \newcommand{\De}{\Delta}

\newcommand{\la} {\lambda}

       \newcommand{\Om}{\Omega}

\newcommand{\Div}{\mathrm{div}}

\begin{document}

\title[Isoperimetric Inequalities and Eigenvalue Estimates]
{Some Isoperimetric Inequalities and Eigenvalue Estimates in Weighted Manifolds}

\author{M. Batista}     
\address{IM, Universidade Fe\-deral de Alagoas, Macei\'o, 
AL, CEP 57072-970, Brazil}
\email{mhbs@mat.ufal.br}

\author{M. P. Cavalcante}     
\email{marcos@pos.mat.ufal.br}

\author{J. Pyo}     
\address{Department of Mathematics  Pusan National University, Busan 609-735, Korea} 
\email{jcpyo@pusan.ac.kr}

\date{\today}

\thanks{The first author was supported by FAPEAL/Brazil grant 20110901-011-0025-0044. }
\thanks{The second author was supported by CNPq/Brazil  grants 306131/2012-9 and 475660/2013-7.}
\thanks{The last author was supported by NRF-2010-0022951/Korea.}

\begin{abstract}
In this paper we prove  general inequalities involving
the weighted mean curvature of compact submanifolds 
immersed in weighted manifolds.  As a consequence we 
obtain a relative linear isoperimetric inequality for such submanifolds.
We also prove an extrinsic upper bound to the first non zero eigenvalue of the 
drift Laplacian on closed submanifolds of  weighted manifolds. 
\end{abstract}

\maketitle

\section{Introduction}

%
%

Let $(\bar{M}^d,\bar g,d\bar\mu)$ be a  weighted manifold, that is, a Riemannian manifold  
$(\bar M^d,\bar g)$  endowed with
a weighted volume form  $ d\bar\mu=e^{-f}d\bar M$, where $f$ is a real-valued smooth 
function on $\bar{M}$ and $d\bar M$ is the volume element induced by the metric $\bar g$.  

\smallskip
In weighted manifolds a natural generalization of the Ricci tensor is the  $m$-Bakry-\'Emery tensor defined by 
$$
\bar {Ric}_f^m= \bar {Ric} + \bar\nabla^2f - \dfrac{1}{m-d}\, df\otimes df,
$$ 
for each  $m\in [d,\infty).$ When $m=\infty$ it gives the tensor $Ric_f=Ric+\bar \nabla^2 f$ introduced by Lichnerowicz
\cite{l1, l2} and independently  by Bakry and \'Emery in \cite{be}. The case $m=d$ only makes sense when the function $f$ is
constant and so  $\bar Ric_f^m$ is the usual Ricci tensor $\bar Ric$ of $\bar M$.



\smallskip

In this paper we are interested in studying inequalities on submanifolds of weighted manifolds. 
In order to do it we make use of intrinsic objects, like the $m$-Bakry-\'Emery tensor, and extrinsic objets
like the weighted mean curvature defined below.
Namely,  given $x:M\to \bar M$  an isometric immersion,  we define the weighted mean curvature vector
${\bf H}_f$ by 
$$
{\bf H}_f= {\bf H} +\bar{\nabla}f^\perp,
$$ 
where ${\bf H} $ is the  mean curvature vector of the submanifold $M$ and ${}^\perp$ denotes the orthogonal projection 
onto the normal bundle $TM^\perp$ (see Gromov \cite{g} and Morgan \cite{m}). 
The  weighted mean curvature appears naturally in the first 
variation of the weighted area functional as described in  \cite{b}. 
In the submanifold $M$ we also consider the weighted volume given $d\mu = e^{-f}dM$, where $dM$ is the
volume element of $M$. 

\smallskip

In case that  $\bar M = \Omega^{n+1}$, where $\Om$ is  a compact oriented 
$(n+1)$-dimensional Riemannian manifold with smooth boundary   $M^n=\partial\Omega$ we consider
on $M$  the Riemannian metric induced by the inclusion map  $\iota: M\hookrightarrow\Omega.$

Let $\nu$ be a unit normal vector field on $M$ and let  $A$ denote the shape operator 
of $M$, that is $A = -\nabla_{(.)}\nu$. It is easy to see that  ${\bf H}_f=H_f\nu$, where 
$H_f=H+\langle\bar\nabla f,\nu\rangle$ and $H= \textrm{trace}\, A$ is the mean curvature
function.



\smallskip

In \cite{ros1}, Ros proved an inequality relating the volume of $\Om$ and the mean curvature function $H$ of 
$M$. 
The inequality  obtained by Ros is essentially contained in the paper of Heintze and Karcher \cite{hk},
although the proof uses different techniques.  

\smallskip
Our first result is the natural generalization of Ros inequality in the context of weighted manifolds.

\begin{theorem}\label{t1} Let $\Om^{n+1}$ be a compact weighted  manifold with 
smooth boundary $M$ and non-negative $m$-Bakry-Emery tensor. Let $H_f$ be the weighted
mean curvature of $M$. If  $H_f$ is positive everywhere, then 
$$
{Vol}_f(\Omega)\leq \dfrac{m-1}{m}\int_M \dfrac{1}{H_f} d\mu.
$$ 

Moreover,  equality holds if and only if  $\Omega$ is isometric to a Euclidean ball,  
$f$ is constant and $m=n+1.$
\end{theorem}


Extending the Ros formula, Choe and Park \cite{cp} proved that a compact 
connected embedded CMC hypersurface  in a convex Euclidean solid cone which 
is perpendicular to the boundary of the cone is part of a round sphere.

The rigidity of compact submanifold with free boundary is a very classical problem 
in  submanifold theory. For instance,  
Nitsche \cite{n} proved that an immersed disk type constant mean curvature surface in a ball which makes a 
constant angle with the boundary of the ball is part of a round sphere.

On weighted manifolds, Ca\~nete and Rosales  \cite{cr} showed the rigidity of compact \emph{stable} 
hypersurfaces  with free boundary in a convex solid cone in Euclidean space with homogeneous density.   
Our next result extends Choe and Park's result to weighted Euclidean spaces $(\Bbb R^{n+1}, ds_0, d\bar{\mu})$, where $ds_0$ is the Euclidean metric. 

\medskip

\begin{theorem}\label{t2}
Let $C$ be a convex solid cone with piecewise smooth boundary $\partial C$ in a weighted manifold $(\Bbb R^{n+1}, ds_0, d\bar{\mu})$ of non-negative m-Bakry-\'Emery tensor. Let $M$ be a compact connected embedded hypersurface in  $C$  and $\Omega$ the bounded domain enclosed by $M$ and $\partial C$. If the weighted mean curvature $H_{f}$ of $M$ is positive everywhere, then
$$Vol_{f}(\Omega)\leq\frac{m-1}{m}\int_{M}\frac{1}{H_{f}}d\mu.$$

Moreover, the equality holds if and only if $M$ is part of a round sphere centered at the vertex of $C$ and $f$ is constant and $m=n+1$.  
\end{theorem}

\medskip

When the weighted mean curvature $H_{f}$ is constant on $M$, we obtain the following relative linear isoperimetric inequality:

\medskip

\begin{corollary}\label{c1}
In Theorem \ref{t1} or \ref{t2}, if $H_{f}$   a positive constant, 
then $$H_{f}Vol_{f}(\Omega)\leq\frac{m-1}{m}Vol_{f}(M).$$
Moreover, the equality holds if and only if $M$ is part of round sphere with  $f\equiv const.$ and $m=n+1$. 
\end{corollary}


\smallskip
\begin{remark} We point out that Morgan \cite{m} and  Bayle  \cite{b} found others generalizations of
the Heintze-Karcher inequality in the context of weighted manifolds. 
Recently,  Huang and Ruan \cite{hr}, give slightly different proofs of Theorem \ref{t1} and its
Corollary. 


\end{remark}

In the second part of this paper, motivated by the work of Heintze \cite{h} 
we  consider the problem to determine extrinsic upper bounds of the 
first eigenvalue of the $f$-Laplacian on closed submanifolds when  the ambient space has radial 
sectional curvature bounded from above.   
We recall that the $f$-Laplacian on $ M$ is defined by
$\De_f u =\De u - \langle  \nabla f, \nabla u \rangle$ for $u\in H^2(M)$. 
When $M$ is closed, it is a basic fact that the spectrum of $\De_f$ is discrete and 
its first non zero eigenvalue is given by  
$$
\la_1(\De_f) = \inf \bigg\{ \frac{\int_M |\nabla \phi |^2 d\mu}{\int_M \phi^2 \, d\mu}: 
\int_M \phi \, d\mu = 0, \, \phi\in C^\infty (M)  \bigg\}.
$$

Using the above notation we have the following results.

\begin{theorem}\label{t4}
Let $ \bar M^{n+p}$ be a weighted manifold with $\sect\leq \delta$, $\delta < 0$. 
If $x: M^n \to \bar M^{n+p}$ is an isometric immersion of a closed manifold, then
$$
\lambda_1(\Delta_f)\leq n\delta +\dfrac{1}{n} \max_M |{\bf H}_f-\bar\nabla f|^2.
$$
\end{theorem}

\medskip

\begin{theorem}\label{t5}
Let $ \bar M^{n+p}$ be a weighted manifold with $\sect\leq \delta$, $\delta\geq 0$. 
If $x: M^n \to \bar M^{n+p}$ is an isometric immersion of a closed manifold such that $x(M)$ 
is contained in a geodesic ball of radius less or equal to $\dfrac{\pi}{4\sqrt\delta}$, then
$$
\lambda_1(\Delta_f)\leq n\delta + \dfrac{1}{nVol_f(M)}\int_M|{\bf H}_f-\bar\nabla f|^2d\mu.
$$ 
Moreover, if equality holds then $M$ is $f$-minimally immersed  into  
$F^{-1}(c)$, where $F=\la f+ \int^r s_\de(t)\, dt$, for some constants $\la$ and $c$, provided that $c$ is a regular value of $F$.
\end{theorem}

\smallskip
\begin{remark}  
Let $(\mathbb{Q}^{n+1}_\de, g_{can})$ be the space form of curvature $\de$  and let $r_0$ be the positive solution of the 
equation $r s_\de ( r )=2nc_\de( r )$ (see Section \ref{extrinsic} for notations). 
A straightforward computation shows that the geodesic sphere of radius $r_0$ has zero weighted mean curvature
in the weighted manifold  $(\mathbb{Q}^{n+1}_\de, g_{can}, e^{-r^2/4})$ and thus the term $|{\bf H}_f-\bar\nabla f|^2$ in 
Theorems  \ref{t4} and \ref{t5} cannot {be}  replaced by $|{\bf H}_f|^2$ when $\de\leq 0$.

%

\end{remark}

\section{Isoperimetric Inequalities}

In this section we recall some {well-known} results on weighted manifolds and we prove Theorem \ref{t1}
and Theorem \ref{t2}.
The first tool we need  is the following Reilly formula (see \cite{md}).

\begin{theorem}\label{reilly1}
Let $u$ be a smooth function on $\Om$. Then we have
$$
\int_{\Omega}\left((\bar\Delta_fu)^2-|\bar\nabla^2u|^2-\bar Ric_f(\bar\nabla u,\bar\nabla u)\right)d\bar\mu=
\int_M\left(2u_\nu\Delta_fu+u_\nu^2H_f+\langle A\nabla u,\nabla u\rangle\right)d\mu,
$$ 
where $\nu$ is the outward unit normal to $M$.
\end{theorem}



Using the Cauchy-Schwarz inequality we have
\begin{proposition}\label{aux}
Let $u$ be a smooth function on $\Om$. Then we have 
$$|\bar\nabla^2u|^2+\bar Ric_f(\bar\nabla u,\bar\nabla u)
\geq\dfrac{(\bar\Delta_f u)^2}{m}+\bar Ric_f^m(\bar\nabla u,\bar\nabla u),
$$ 
for every $m> n+1$ or $m=n+1$ and $f$ is a constant.  Moreover,  equality holds if  and only if $\bar\nabla^2u=\lambda g$ and 
$\langle\bar\nabla u,\bar\nabla f\rangle=-\dfrac{m-n-1}{m}\bar\Delta_f u.$
\end{proposition}

%
%


\medskip
The last tool is an important result due to Reilly (see \cite{r}).
\begin{theorem}\label{reilly2}
Suppose that $\Omega$ admits a function $u:\Omega\to\real$ and non-zero constant $\lambda$ such that 
\begin{enumerate}
\item[(a)] $\bar\nabla^2u=\lambda g;$
\item[(b)] $u|_M$ is constant.
\end{enumerate}
Then $\Omega$ is isometric to an Euclidean ball.
\end{theorem}

\subsection{Proof of Theorem \ref{t1}}
Let $u:\Omega\to\real$ be the solution of the Dirichlet problem \\
\begin{equation}\label{sistema}
\left\{ \begin{array}{rcc}
\bar\Delta_f u=1& \mbox{in} & \Omega, \\ 
u=0 & \mbox{in} & \partial\Omega . 
\end{array}\right.
\end{equation}
\\
Plugging this function in Theorem \ref{reilly1}  we get 
$$
\int_{\Omega}\left(1-|\bar\nabla^2u|^2-\bar Ric_f(\bar\nabla u,\bar\nabla u)\right)d\bar\mu=\int_Mu_\nu^2H_f \, d\mu.
$$ 
Using the Proposition \ref{aux} we have
$$
\int_{\Omega}\left(1-\dfrac{1}{m}-\bar Ric^m_f(\bar\nabla u,\bar\nabla u)\right)d\bar\mu\geq\int_Mu_\nu^2H_f \, d\mu.
$$ 
Using the hypothesis on the $m$-Bakry-\'Emery tensor we obtain
$$
\int_Mu_\nu^2H_f \, d\mu\leq \dfrac{m-1}{m}Vol_f(\Omega).
$$
\\
Hence, using the Stokes theorem and the above inequality we have

\begin{equation*} 
\begin{array}{lll}
Vol_f(\Omega)^2&=&\displaystyle \medskip\left(\int_\Omega\bar\Delta_fu\, d\bar\mu\right)^2 
=  \left(\int_Mu_\nu\, d\mu\right)^2 \\
&=& \displaystyle \medskip  \left(\int_Mu_\nu\sqrt{H_f}(\sqrt{H_f})^{-1}\, d\mu\right)^2\\
& \leq& \displaystyle \medskip \left(\int_Mu_\nu^2H_f\, d\mu\right)\left(\int_M\dfrac{1}{H_f}\, d\mu\right)
\leq\dfrac{m-1}{m}Vol_f(\Omega) \int_M\dfrac{1}{H_f}\, d\mu.
\end{array} 
\end{equation*}
\\
That is,
\begin{equation}\label{kh}
Vol_f(\Omega)\leq\dfrac{m-1}{m} \int_M\dfrac{1}{H_f}\, d\mu .
\end{equation}

Now assume that equality occurs in (\ref{kh}). Then all the inequalities above are equalities and thus we obtain  
\begin{equation}\label{id}
\left\{  \begin{array}{lll}
\medskip \bar\nabla^2u=\lambda g, \\ 
\medskip
\langle\bar\nabla u,\bar\nabla f\rangle=-\dfrac{m-n-1}{m}\bar\Delta_f u, \\ 
\bar Ric_f^m(\bar\nabla u,\bar\nabla u)=0. 
   \end{array} \right.
\end{equation}

From the first   equation above, we have $\bar \De u= (n+1)\lambda$. So, using the
first equation in (\ref{sistema}) and the second equation above it is easy to
see that $\lambda = \frac{1}{m}$. 
Since $\la$ is constant we apply Theorem \ref{reilly2} to obtain the $\Om$ is isometric to a
Euclidean ball. 

\smallskip
Finally, assume that $m>n+1$. Then, the first equation in (\ref{id})  and the 
boundary condition in (\ref{sistema})  imply that  $u=\frac{\lambda}{2}r^2+C$, for some constant $C$,
where $r$ is the distance function from its minimal point (see \cite{r}). 
Therefore, using the second equation in (\ref{id}) we find $f=-(m-n-1)\ln r +C$.
It is a contradiction, since $f$ is a smooth function.

\subsection{Proof of Theorem \ref{t2}}

Let $\Omega_{\epsilon}\subset\Omega$ be a domain with smooth boundary which is 
obtained from $\Omega$ 
by smoothing out the region within a distance  $\epsilon>0$ from
the singular set of $\partial \Omega$. Let $u$ be a smooth solution to the following mixed 
boundary value problem:

\begin{equation}\label{mixed}
\left\{  \begin{array}{lll}
\medskip \bar{\triangle}_{f} u=1 & \textrm{in $\Omega_{\epsilon}$}\\
\medskip
u=0 & \textrm{on $\partial\Omega_{\epsilon}\setminus \partial C$}\\ 
u_{\nu}=0  & \textrm{on $\partial\Omega_{\epsilon}\cap \partial C$}.   \end{array} \right.
\end{equation}

Applying the solution of the above problem $u$  into Theorem 2.1 and using  Proposition \ref{aux} we get
$$
\int_{\Omega_{\epsilon}}
\left(1-\dfrac{1}{m}-\bar Ric^m_f(\bar\nabla u,\bar\nabla u)\right)d\bar\mu
\geq\int_{\partial\Omega_{\epsilon}\setminus\partial C}u_{\nu}^2H_{f}d\mu+
\int_{\partial\Omega_{\epsilon}\cap\partial C}\langle A\nabla u, \nabla u\rangle d\mu.$$ 

Since $C$ is convex, $\int_{\partial\Omega_{\epsilon}\cap\partial C}\langle A\nabla u, \nabla u\rangle d\mu\geq0$.
Now, using that $\bar Ric ^m_f \geq 0$ we obtain

$$
\frac{m-1}{m}Vol_f(\Om_\epsilon)
\geq\int_{\partial\Omega_\epsilon\setminus\partial C}u_{\nu}^2H_{f}d\mu.$$

By the same argument as in the proof of Theorem \ref{t1}, we get 
\begin{equation*}
Vol_f(\Omega_\epsilon)\leq\dfrac{m-1}{m} \int_{\partial\Omega_\epsilon\setminus\partial C}\dfrac{1}{H_f}\, d\mu .
\end{equation*}

Letting $\epsilon \to 0$ we obtain the desired inequality. Here, it is important to point out that
$H_f\to \infty$ near $\partial M \cup \hat C$ since
$H\to \infty$ and $\bar \nabla f $ is bounded in $\Om$, where $\hat C$ is the singular set of $C$,  (see \cite{cp}).

Now we assume that equality holds, then we get $f$ is constant and $m=n+1$ by the same argument as 
in the proof of  Theorem \ref{t1}.
Let $O$ be the vertex of the solid cone $C$.
For a constant $R>0$, $u(X)=\frac{|X-O|^2-R^2}{2(n+1)}$ is the solution of the mixed boundary value problem (\ref{mixed}).
Since $u(X)=0$ on $M$, $M$ is part of a round sphere centered at $O$ and $M$ meets $\partial C$ with right angle along the boundary.

By simple computation, the converse holds. This completes the proof of Theorem \ref{t2}.

\section{Extrinsic Eigenvalue Estimates}\label{extrinsic}


Let $x: M^n \to \bar M$ be an isometric immersion of a closed manifold.
Given $Y$ a vector field on $\bar M$, we denote by $D_f Y $  the (extrinsic) $f$-divergence
of $Y$, that is
$$
D_fY := \Div_M Y - \langle \bar \nabla f, Y \rangle.
$$
Note that if $u$ is a smooth function on $M$, then $D_f (\nabla u) = \De_f u$.
%

\smallskip

In the sequel we assume the radial sectional curvature of  $\bar M$ is bounded from above, that is,
there exists a constant $\de$ such that $\sect\leq \de$.  
Let us consider the vector field $X=s_{\de}(r)\bar\nabla r$ on $\bar M$, where the function
$s_{\de}$ is the solution of the ODE
\begin{equation*}
\left\{ \begin{array}{lll}
 \medskip  g''(t)+\de g(t)=0  \\
g(0)=0, \quad  g'(0)=1  .
\end{array}\right.
\end{equation*}

In the lemma below $c_\de$ denotes the derivative of $s_\de$ and
$X^\top$ is the tangent component of $X$ on $M$.

\begin{lemma}\label{lemma1} 
 On the above conditions  we have:
\begin{enumerate}
\medskip
\item  $D_fX  \geq n c_\de - s_\de\langle \bar \nabla f, \bar\nabla r  \rangle;$
\medskip
\item $D_f X^\top \geq nc_\de - s_\de \langle \mathbf{ H_{f}} -\bar \nabla f, \bar \nabla r \rangle;$
\medskip
\item $n\int_M c_\de  d\mu \leq - \int_M s_\de \langle \mathbf{ H_{f}}-\bar \nabla f, \bar \nabla r \rangle d\mu
\leq
\int_M s_\de  |\mathbf{ H_{f}}-\bar \nabla f| d\mu.$ 
\end{enumerate}
\end{lemma}
The proof is a slight modification of Lemmas 2.4 and 2.5 of Heintze \cite{h},  using the 
weighted volume in assertion (3). We point out that if equalities hold in assertion (3),
then there is a function $\la$ on $M$ such that 
$X(p )=\la(p ) (\mathbf{ H_{f}}-\bar \nabla f)(p )$, $\forall p\in M.$

\begin{lemma}\label{lemma2}  On the above conditions  we have:
$$
\de \int_M |X^\top|^2 d\mu \geq n \int_M c_\de^2\, d\mu - \int_M c_\de \, s_\de\, |\mathbf{ H_{f}} -\bar \nabla f |\, d\mu.
$$
\end{lemma}

{\it Proof.}
The case  $\de=0$ follows from assertion (3) in the previous  lemma. 

\smallskip
Assume $\de\neq 0$ and note that $X^\top = \nabla(-\frac{c_\de}{\de})$. Using Lemma \ref{lemma1} we have
\begin{equation*}
\begin{array}{lll} \medskip
\displaystyle\de\int_M |X^\top|^2 d\mu &=&\displaystyle \de\int_M |\nabla(-\frac{c_\de}{\de})|^2 d\mu \\ \medskip
&=& \displaystyle\de\int_M \frac{1}{\de} c_\de \, \De_f(-\frac{c_\de}{\de}) d\mu \\ \medskip
&=& \displaystyle \int_M  c_\de \, D_f(s_\de \nabla r) d\mu \\ \medskip
&=&  \displaystyle\int_M  c_\de \, D_f \, X^\top d\mu \\  \medskip
&\geq& \displaystyle \int_M  c_\de \big( nc_\de - s_\de \langle \mathbf{ H_{f}}-\bar \nabla f, \bar \nabla r \rangle \big) d\mu \\ 
&\geq& \displaystyle n\int_M  c_\de^2 d\mu - \int_M  c_\de  s_\de | \mathbf{ H_{f}}-\bar \nabla f | d\mu.
\end{array}
\end{equation*}


\bigskip

The next lemma is the compilation of Lemmas 2.7  and 2.8 of \cite{h}. 

\begin{lemma}\label{lemma3} On the above notations the following assertions hold:
\begin{enumerate}
\item 
Let $(x_1, \ldots, x_m)$ denote a system of  normal coordinates on $\bar M$ around the 
center of mass of $M$ in $\bar M$ with respect to the weighted volume.   Then we have
$$
\sum_{i=1}^m \big|\nabla \big( \frac{s_\de}{r} x_i\big)\big |^2 + \de|X^\top|^2 \leq n.
$$

\item If $\de \leq 0$,
$$
\int_M s_\de d\mu \, \int_M s_\de\, c_\de d\mu \leq \int_M s_\de^2 d\mu \, \int_M c_\de d\mu.
$$
\end{enumerate}
\end{lemma}

\medskip
Notice that for each $i=1,\ldots, n$, we have 
$
\int_M  \frac{s_\de}{r} x_i \, d\mu = 0.
$

\subsection{Proof  of Theorem  \ref{t4} } 
Firstly, we write  $s_\de^2 = \sum \frac{s_\de^2}{r^2} x_i^2 $. Using
the functions $ \frac{s_\de}{r} x_i$ as test  functions  in 
the variational characterization of $\la_1(\De_f)$
and Lemmas \ref{lemma2} and \ref{lemma3} we have
\begin{equation*}
\begin{array}{lll}
\medskip
\displaystyle
\la_1(\De_f)\, \int_M s_\de^2 d\, \mu  &\leq&\displaystyle \smallskip\int_M \sum \big|\nabla  \frac{s_\de }{r} x_i\big|^2\, d\mu \\
&\leq&\displaystyle \smallskip  \int_M (n - \de |X^\top|^2)\, d\mu \\
&\leq& \displaystyle \smallskip n Vol_f (M) - n \int_M c_\de^2\, d\mu + \int_M c_\de \, s_\de\, |\mathbf{ H_{f}} -\bar \nabla f |\, d\mu\\
&\leq& \displaystyle \smallskip n  \int_M \de s_\de^2\, d\mu + \max_{M} |\mathbf{ H_{f}} -\bar \nabla f |\, \int_M c_\de \, s_\de\,  d\mu\\
&\leq& \displaystyle \smallskip n  \int_M \de s_\de^2\, d\mu + \max_{M} |\mathbf{ H_{f}} -\bar \nabla f |\, \int_M s^2_\de \,  d\mu
\, \frac{\int_M c_\de \,  d\mu}{\int_M s_\de\,  d\mu}\\
&\leq& \displaystyle \smallskip n  \int_M \de s_\de^2\, d\mu + \max_{M} |\mathbf{ H_{f}} -\bar \nabla f |\, \int_M s^2_\de \,  d\mu
\, \frac{\frac{1}{n}\int_M s_\de\,|\mathbf{ H_{f}}-\bar \nabla f| \,  d\mu}{\int_M s_\de\,  d\mu}\\
&\leq& \displaystyle \smallskip n  \int_M \de s_\de^2\, d\mu + \frac{1}{n}\max_{M} |\mathbf{ H_{f}} -\bar \nabla f |^2\, \int_M s^2_\de \,  d\mu
\end{array}
\end{equation*}
and the result follows.

\bigskip 
\subsection{Proof  of Theorem  \ref{t5}}
We give the proof in two cases.

Case 1: $\de=0$. Take the functions $x_i$ given in the Lemma \ref{lemma3} as test function in the
 variational characterization of $\la_1(\De_f)$. Taking the sum we have,
\begin{equation*}
\la_1(\De_f)\, \int_M  \sum_{i} x_i^2 d\mu  \leq \int_M \sum_{i}|\nabla x_i|^2 d\mu.
\end{equation*}
 
From assertion (1) of  Lemma \ref{lemma3} we have $\sum_{i}|\nabla x_i|^2 \leq n$. So, using
the assertion (3) in
Lemma \ref{lemma1}  we obtain

\begin{equation*}
\begin{array}{lll}
\medskip
\displaystyle
\la_1(\De_f)\, \int_M | X|^2 d\mu  &\leq&  nVol_f(M) \\ 
&\leq& \displaystyle \smallskip\int_M r  |\mathbf{ H_{f}}-\bar \nabla f|\, d\mu \\
&\leq& \displaystyle \smallskip \frac{1}{nVol_f(M)}\bigg(\int_M |X|  |\mathbf{ H_{f}}-\bar \nabla f| \, d\mu\bigg)^2 \\
&\leq& \displaystyle \frac{1}{nVol_f(M)} \int_M |X|^2 d\mu  \int_M |\mathbf{ H_{f}}-\bar \nabla f|^2\, d\mu. 
\end{array}
\end{equation*}

That is,
$$
\la_1(\De_f) \leq \displaystyle \frac{1}{nVol_f(M)}  \int_M |\mathbf{ H_{f}}-\bar \nabla f|^2\, d\mu. 
$$

\bigskip

Case 2: $\de >0.$ In this case let us use as test functions the functions $\frac{s_\de}{r}x_ i$, $i=1,\ldots, m$ and 
the function $\frac{c_\de-c}{\sqrt \de}$, where $c= \frac{1}{Vol_f(M)}\int_M c_\de\, d\mu$. Applying these functions
to the variational characterization of $\la_1(\De_f)$ and using the assertion (1) in Lemma \ref{lemma3} we have

\begin{equation*}
\begin{array}{lll}
\medskip
\displaystyle
\la_1(\De_f)\, \int_M  \big(s_\de^2+ \frac{(c_\de-c)^2}{\de} \big) d\mu  &\leq& 
\displaystyle \smallskip \int_M  \sum_i \big| \nabla \frac{s_\de}{r}x_i \big|^2  +\frac{1}{\de}|\nabla c_\de|^2\, d\mu\\
&=&\displaystyle \smallskip \int_M  \sum_i \big| \nabla \frac{s_\de}{r}x_i \big|^2  +{\de}|X^\top|^2\, d\mu\\
&\leq& nVol_f(M). 

\end{array}
\end{equation*}

On the other hand, from a direct computation we get
$$
 \int_M  \big(s_\de^2+ \frac{(c_\de-c)^2}{\de} \big) d\mu = \dfrac{1-c^2}{\de}Vol_f(M).
$$
So
$$
\la_1(\De_f)(1-c^2)\leq n\de.
$$

\medskip
To finish the proof, we will estimate the term $1-c^2$ from below. 
We set 
$$
d=1+\frac{n^{-2}}{\de  Vol_f(M)}\int_M|H_f-\bar\nabla f|^2\, d\mu.
$$ 
Then we use the assertion (3) in  Lemma \ref{lemma1}  to get

\begin{equation*}
\begin{array}{lll}
\medskip
\displaystyle
 
(1-c^2)d & = & d -  \displaystyle\frac{1}{Vol_f(M)^2}\bigg(\int_M c_\de\, d\mu\bigg)^2\\
&& - \displaystyle \frac{n^{-2}}{\de  Vol_f(M)^3}\bigg(\int_M c_\de\, d\mu\bigg)^2 \int_M|H_f-\bar\nabla f|^2\, d\mu   \\
%
&\geq& 
\smallskip \displaystyle 
d
-  \frac{n^{-2}}{ Vol_f(M)^2}\bigg(\int_M s_\de |H_f-\bar\nabla f|\, d\mu\bigg)^2  \\
&& - \displaystyle \frac{n^{-2}}{\de  Vol_f(M)^2}\int_M c_\de^2\, d\mu  \int_M|H_f-\bar\nabla f|^2\, d\mu. \\

%
&\geq& \smallskip
\displaystyle 1+ 
{n^{-2}}\int_M|H_f-\bar\nabla f|^2 d\mu\,\bigg(\frac{1}{\de Vol_f(M)}-\frac{1}{\de Vol_f(M)^2}\int_M (\de s_\de^2 + c_\de^2)\, d\mu\bigg)\\
\smallskip&=&1.

\end{array}
\end{equation*}
\medskip
Thus 
\begin{equation}\label{final}
\la_1(\De_f)
\leq n\de + \frac{1}{ n Vol_f(M)}\int_M|H_f-\bar\nabla f|^2\, d\mu
\end{equation}
as we claimed. 

\medskip

In what follows, we analyze the case of equality in (\ref{final}). 
In this case all the inequalities above are equalities. So, when the equalities hold
we get
$$
s_\de\bar\nabla r=\la (H_f-\bar\nabla f),
$$
for some constant $\la$. 
In particular,
$$
s_\de\nabla r = -\la \nabla f.
$$
It means that  the function $F=\la f + \int^r s_\de(t)\, dt$ is  constant on $M$. 
Thus $M$ is immersed in $M_0=F^{-1}(c)$ for some constant $c$. If $c$ is a regular value of $F$, 
then we can decompose the second fundamental form of the immersion $x$ as
$$
\alpha = \alpha_0 + \tau,
$$
where $\al_0$ is the  second fundamental form of $M$ in $M_0$ and
$\tau$ is parallel to $\bar\nabla F$.

Let $H_0$ be the mean curvature vector of $M$ in $M_0$. Then
$$
H_0=H- trace \, \tau. 
$$
Denote by $ \widehat \nabla f$ the gradient of $f$ in $M_0$ and by $\widehat\nabla f^\perp$ its normal component on $M$.
Thus,
$$
H_{0}+ \widehat\nabla f^\perp = H_f - \frac{\langle \bar\nabla f,\bar\nabla F\rangle}{|\bar \nabla F|^2}\bar\nabla F- trace(\tau).
$$
It is easy to see that, on $M$, the right hand side of the equality above is parallel to $\bar\nabla F.$ 
Since $H_{0}+ \widehat\nabla f^\perp $ is tangent to $M_0$ 
we conclude that  it vanishes, as we claimed.

\section*{Acknowledgments} 
The authors would like to thank the referee for the very valuable comments and detailed corrections.

\end{document}